\documentclass[english,parskip=half]{scrartcl} 

\usepackage[utf8]{inputenc} 

\usepackage{babel} 
\usepackage{lmodern} 
\usepackage[T1]{fontenc} 
\usepackage{microtype} 
\usepackage{graphicx} 
\usepackage{csquotes} 
\usepackage{verbatim} 
\usepackage[output-decimal-marker={,},exponent-product=\cdot]{siunitx} 
\usepackage[sorting=nyt, giveninits=true, backref=false, isbn=false, doi=false, eprint=false, url=false]{biblatex} 
\renewbibmacro{in:}{}
\usepackage{tabularx} 
\usepackage{url}
\usepackage{mathtools,amssymb} 
\usepackage{amsmath}   
\usepackage{amsthm}    
\usepackage{bbm}
\usepackage{braket}
\usepackage{enumitem}
\usepackage{relsize}
\usepackage[section]{placeins} 
\usepackage{listings}
\usepackage{tikz}
\usepackage{pgfplots}
\usepackage{pgf}
\usepackage{pgfplotstable, booktabs} 
\usepackage{hhline}
\usepackage{xcolor,colortbl}
\usepackage{longtable}
\usetikzlibrary{hobby}
\usepackage{enumitem}
\usepackage{soul} 
\usepackage{algorithm}
\usepackage{algorithmic}

\usepackage{xfrac}

\usepackage{afterpage}

\usepackage{accents}

\def\XXint#1#2#3{{\setbox0=\hbox{$#1{#2#3}{\int}$ }
		\vcenter{\hbox{$#2#3$ }}\kern-.6\wd0}}

\DeclarePairedDelimiter{\abs}{\lvert}{\rvert}
\DeclarePairedDelimiter{\norm}{\lVert}{\rVert}
\NewDocumentCommand{\normL}{ s O{} m }{%
	\IfBooleanTF{#1}{\norm*{#3}}{\norm[#2]{#3}}_{L^2(\Omega)}%
}

\DeclareFontFamily{U}{mathc}{}
\DeclareFontShape{U}{mathc}{m}{it}%
{<->s*[1.03] mathc10}{}

\DeclareMathAlphabet{\mathscr}{U}{mathc}{m}{it}
\usepackage{hyperref}
\usepackage{cleveref}

\newtheorem{theorem}{Theorem}[section]
\newtheorem{remark}[theorem]{Remark}
\newtheorem{assumptions}[theorem]{Assumptions}

\newtheorem{lemma}[theorem]{Lemma}

\crefname{remark}{Remark}{Remarks}
\crefname{assumptions}{Assumptions}{Assumptions}

\bibliography{paper-bib.bib} 

\addto\extrasenglish{%
} 

\makeatletter	
\g@addto@macro\@floatboxreset\centering
\makeatother

\definecolor{whitesmoke}{rgb}{0.96, 0.96, 0.96}
\definecolor{timberwolf}{rgb}{0.86, 0.84, 0.82}

\def\johannes#1{#1}
\def\kommentar#1{}
\def\revone#1{#1}
\def\revtwo#1{#1}

\begin{document}


%
%
%
%
%


\begin{center}
	{\Large\bfseries Optimal control of a quasilinear parabolic equation\\and its time discretization}
	
	Luise Blank\footnote{
		Department of Mathematics,
		University of Regensburg, D-93040 Regensburg, Germany (luise.blank@ur.de, johannes.meisinger@ur.de)}, 
	Johannes Meisinger\footnotemark[1]
\end{center}

\bigskip

\begin{abstract}
\textbf{Abstract.}
In this paper we discuss the optimal control of a quasilinear
parabolic state equation.
{Its form is leaned on the} kind of problems arising for example when controlling the anisotropic Allen-Cahn equation as a model for crystal growth.
Motivated by this application we consider the state equation as a result of a gradient flow of an energy functional. The quasilinear term is strongly monotone and obeys a certain growth condition
and the lower order term is non-monotone.
The state equation is discretized implicitly in time with piecewise constant functions.
The existence of the control-to-state operator and its Lipschitz-continuity is shown for the time discretized as well as for the time continuous problem.
Latter is based on the convergence proof of the discretized solutions.
Finally
we present for both the existence of global minimizers.
Also convergence of 
a subsequence of time discrete optimal controls to a
global minimizer of the time continuous 
problem can be shown.
Our results hold in arbitrary space dimensions.
\end{abstract}
\bigskip

\noindent \textbf{Key words. }
	quasilinear parabolic equation, Allen-Cahn equation, anisotropy, optimal control, implicit discretization, convergence analysis

\noindent \textbf{AMS subject classification. } 35K59, 49J20, 49M41, 65M12, 65M60

\bigskip


\section{The optimization problem}
\label{sec:intro}
In many areas 
the optimal control of an interface evolution
towards an anisotropic shape is desired.
For example in {chemistry or} materials science one wishes to steer the solidification process of crystals \cite{BorNRVS09,EisVS14,fujiwara2006growth,nakajima2009crystal}.
For the time evolution of shapes phase-field models have shown great promise in many application areas and anisotropies can be incorporated
{(see, e.g., \cite{DDE} and references therein).}
In this ansatz the interface is modeled with a  diffuse interface layer,
and an order parameter $y$---the so called phase-field---%
reflects the pure phases with the values
$\pm1 $, e.g. the liquid phase for
$y \approx 1$ and the solid phase when $y \approx -1$, and the diffuse interface with values between $-1$ and $1$.
%
The gradient flow of 
the Ginzburg-Landau energy, having the form
\begin{equation}
\label{GinzburgLandauEnergyAniso}
\mathcal{E}(y) \coloneqq \int_{\Omega} 
A(\nabla y) + 
\psi(y) \, \mathrm dx
\end{equation}
then determines the time evolution of the shape, 
and with it the state equation for the control problem.
Here the first term represents the surface energy where
$A:\mathbb{R}^d \to \mathbb{R}$ is an (an-)isotropy function,
and the potential
$\psi$ 
can be thought of being symmetric and to have its global minima at $\approx \pm 1$.
Let us mention that typically the energy involves a variable $\varepsilon >0$
related to the interfacial thickness which we set to 1 without loss of generality in this paper.
Considering in particular the $L^2$-gradient flow with a smooth potential we obtain the Allen-Cahn equation.
For further introduction to phase field models we refer to {\cite
{eck2017mathematical}} and references therein.
The following analysis and numerical ansatz will not only be valid for the Allen-Cahn equations but can be applied  in general to differential equations arising from a gradient flow of energies of the form~\cref{GinzburgLandauEnergyAniso}.
\revone{For the gradient of $A$ we use the notation $A'(p)$ instead of $\nabla A(p)$ as is common for anisotropic phase field models.}

The goal is now to determine the distributed control $u$ driving the solution $y$ of the
gradient flow equation
\begin{equation}
  \partial_t y - \nabla \cdot A'(\nabla y)  +
  \psi'(y)  =  u 
\end{equation}
from an initial configuration $y_0$ at time $t_0$---say $t_0=0$---to a given target function $y_\Omega$ at a given final time $T$
(or a target function  $y_Q$ over the whole time horizon as considered in \cref{eq:distributed_cost_disc}).
Hence the optimal control problem is described by the following setting:
\\
Let $\Omega \subset \mathbb{R}^d$ 
be a bounded Lipschitz domain and $y_\Omega \in L^2(\Omega)$ be a given target function.
Let a final time $0<T<\infty$ be given and denote the space-time cylinder by $Q \coloneqq [0,T]\times \Omega$ and its \revone{lateral} boundary by $\Sigma \coloneqq [0,T] \times \partial \Omega$.
Our objective is to find for a given initial state
$y_0 \in H^1(\Omega)$
a solution to the optimal control problem
\begin{equation}
\label{eq1a}
\min J(y,u) \coloneqq \dfrac{1}{2} \norm{y(T) - y_\Omega}^2_{L^2(\Omega)} + \dfrac{\lambda}{2} \norm{u}^2_{L^2(Q)},
\end{equation}
subject to the quasilinear, possibly nonsmooth parabolic state equation
\begin{equation}\label{weakAC}
\begin{aligned}
  \int_Q \partial_t y \eta + A'(\nabla y)^T \nabla \eta +
  \psi'(y) \eta &=  \int_Q u\eta 
     \quad  \forall \eta \in L^2(0, T; H^1(\Omega)),\\
  y(0)&=y_0     \quad \mbox{ in } \Omega,
\end{aligned}
\end{equation}
where $u\in L^2(Q)\cong L^2(0,T;L^2(\Omega))$ and $y \in L^2(0,T;H^1(\Omega))\cap H^1(0,T;H^1(\Omega)')$.
Note that $J$ is well defined due to the embedding $L^2(0,T;H^1(\Omega))\cap H^1(0,T;H^1(\Omega)') \hookrightarrow C([0,T]; L^2(\Omega))$.
\revone{The weak formulation implies the boundary condition $A'(\nabla y)^T\nu=0$ on $\Sigma$ where $\nu$ is the outer normal.}

In contrast to the well-studied isotropic Allen-Cahn equation where $A'=Id$,
in the anisotropic case
$A:\mathbb{R}^d \to \mathbb{R}$ is an absolutely 2-homogeneous function.
\kommentar{Würde hier keinen Literaturverweis machen. Wenn dann Deckelnick/Dziuk.., Kobayashi, ElliottSchaetzle, BGN3}%
As a consequence $A'$ is not differentiable at $0$ in general.
Let us mention that roughly speaking $y$ is constant in the pure phases, and hence
$\nabla y \approx 0$ holds except on the interface. 
However, for solving an optimal control problem numerically one typically uses the differentiability of the control-to-state operator which would require differentiability of $A'$.
\johannes{This paper serves also as a preparation to \cite{BlMeNumeric}, where we show the differentiability and first order conditions of an implicit in time discretized problem with a regularized $A$.}
\johannes{Hence,} to allow the use of possibly regularized
anisotropies $A$
we relax the requirement of 2-homogeneity.
Moreover we
account for various potentials $\psi$.
\begin{assumptions}\label{assu}
Assume $A\in C^1(\mathbb{R}^d)$
with $A'$ being strongly monotone, i.e.,
\begin{equation*}
(A'(p)-A'(q),p-q) \geq C_A \vert p-q\vert^2 \qquad \forall p,q\in\mathbb{R}^d,
\end{equation*}
{with $C_A>0$} and such that it fulfills the growth condition $ \vert A'(p)\vert \leq \overline{C}_A\vert p \vert $ with $\overline{C}_A>0$.\\
Furthermore let $\psi \in C^1(\mathbb{R})$
be bounded from below and \revone{such that it} can be approximated by 
$f_n $ satisfying 
\begin{equation}
\label{f_approximation_property}
f_n \in C^2(\mathbb{R}),\quad
f_n \to \psi \quad \text{in } C^1_{\text{loc}},\quad
-c \leq f_n \leq c (\psi+1) , \quad 
f_n'' \geq -C_\psi, \quad \vert f_n'' \vert \leq C_n ,
\end{equation}
with $c, C_n,C_\psi \geq 0$ 
and  $\psi(y_0)\in L^1(\Omega)$ for the given initial data $y_0\in H^1(\Omega)$.
\end{assumptions}
\revtwo{In the following we use the convention $\tfrac{1}{C_\psi}\coloneqq\infty$ for $C_\psi=0$.}
\revone{Note that $A'$ being strongly monotone on $\mathbb{R}^d$ is equivalent to $A$ being strongly convex. Further, the}
assumptions on $\psi$ in particular imply that it holds
	\begin{equation}
	\label{eq:lower_estimate_psi_prime_diff}
		\left(\psi'(y_1) - \psi'(y_2), y_1-y_2\right) \geq -C_\psi \vert y_1-y_2\vert^2 \qquad \forall y_1, y_2 \in \mathbb{R}.
	\end{equation}
        \kommentar{
sketch for (*):\\
$"\implies":$ $x = y+tv$
$$\left(\psi'(y+tv) - \psi'(y)\right)tv \geq C t^2v^2 \iff \left(\tfrac{\psi'(y+tv) - \psi'(y)}{t}\right)v \geq C v^2 \stackrel{\lim_{t\to 0}}{\implies} \psi''(y) v^2 \geq Cv^2 \iff \psi''(y) \geq C$$
$"\Longleftarrow":"$ by mean value theorem there exists $\xi_{x,y} \in [x,y]$ with $\psi''(\xi_{x,y}) (x-y) = \psi'(x) - \psi'(y)$.\\
$\implies$ $\left(\psi'(x) - \psi'(y)\right)(x-y) = \psi''(\xi_{x,y}) (x-y)^2 \geq C (x-y	)^2$\\
If $\psi$ is only approximated, we can use the above to show the inequality for the $f_n$. Then for $y_1, y_2$ fixed it holds
\begin{align*}
  &\left(\psi'(y_1) - \psi'(y_2), y_1-y_2\right) \\
  &= \left((\psi'(y_1)-f_n'(y_1) + (f_n'(y_2)-\psi'(y_2)) + (f_n'(y_1) - f_n'(y_2)), y_1-y_2\right) \\
&\geq	\left((\psi'(y_1)-f_n'(y_1) + (f_n'(y_2)-\psi'(y_2)), y_1-y_2\right) - C_\psi |y_1-y_2|^2.
\end{align*}
Taking the limit using $f_n' \to \psi'$ in $C_{\text{loc}}$ then shows \cref{eq:lower_estimate_psi_prime_diff} (Note that $y_1,y_2$ were fixed but the reasoning can be done for all $y_1,y_2$ separately).
}
Some examples of $A$ and $\psi$ with respect to Allen-Cahn equations are mentioned in \cref{possiblepsi}.

In this paper we study the existence of an optimal control to  \cref{eq1a}--\cref{weakAC} in arbitrary space
dimension,
the existence of \johannes{a solution to} the corresponding in time discretized control problem,
and the convergence of the time-discrete optimal controls to \johannes{a} time-continuous one.
Here an implicit time discretization using
piecewise constant functions is employed.
Under some additional smoothness requirements on $A$ and $\psi$ 
the first order optimality condition
is adressed in  \cite{BlMeNumeric}
where its derivation  relies on results of this paper. 
Therein one can \johannes{also} find numerical results.
 
To the
best of our knowledge there does not exist any mathematical treatment on the optimal control of anisotropic phase-field models so far.
%
Optimal control of isotropic Allen-Cahn variational equations are studied, e.g., in \cite{Benner2013,BFSHMR, CFS15, farshbaf2013optimal, OSY12} and of Cahn-Hilliard variational {(in-)}equalities in \cite{HINTERMULLER2013810, CFSGS15, HW2012} and references
therein.
%
Let us mention results given in the context of anisotropic Allen-Cahn equations.
One possible anisotropy was introduced in a pioneering paper by Kobayashi \cite{Kobayashi} and existence and uniqueness of a solution are studied in \cite{Burman-Rappaz, Miranville, TLVW}.
For quite general anisotropies the solution \revtwo{to} Allen-Cahn equations with obstacle potential is analyzed in  \cite{Elliott1996}. Among others they use 2-homogeneity of $A$, an approximation of the potential similar to \cref{f_approximation_property} and an implicit time discretization (without showing convergence of the discretization).
Explicit and semi-implicit approximations
are discussed in the survey paper \cite{DDE}, where also many references are given.
For convex Kobayashi anisotropies
several time discretizations  are considered in \cite{Graser2013}. 
In \cite{Barrett2014, bgn13} particular suggestions for the anisotropies are given and an efficient semi-implicit method using a particular linearization of $A'$ and a convex/concave splitting
is presented and energy stability is shown. Also several numerical experiments are shown comparing the anisotropies. 
\\
Literature to optimal control of quasilinear parabolic equations of the form \cref{weakAC} is still in its infancy.
Most literature known to us treat quasilinearities with coefficients depending on $x,t$ 
{and on the function $y$ but not on its gradient} 
\cite{Neitzel2018, Neitzel2020, Casas2018, Meinlschmidt2017a, Meinlschmidt2017b}.
For quasilinearities involving spatial derivatives of $y$ see for example \cite{TroelSerge2017, Casas1995parabolic}.
In particular let us mention that the latter reference
contains the most similar problem to ours, as the authors require a rather general quasilinearity with some
particular polynomial growth condition.
However
they require the nonlinearity
{$\psi'$} to be monotone.
All the literature listed here assumes the quasilinear term to be rather well behaved, in particular none of its derivatives shall be singular at the origin.
In the present context to our knowledge such difficulties have only been considered for elliptic equations~\cite{Casas1991}.

%

The outline of the paper is the following.
\\
As a first step we study the state equations.
Therefore we introduce in \cref{sec:analysis_continuous}
the time discretization.
Then we discuss the existence and uniqueness
of the solution of the discretized state equation
as well as the Lipschitz-continuity of the control-to-state operator. Furthermore, for a set of bounded controls we obtain bounds on the states \revtwo{independent} of the discretization level.
Using these results we consider the limit with respect to the time discretization and obtain corresponding results for the in time continuous state equation \cref{weakAC}. Consequently we have also convergence of the discretization. In addition we show energy stability of the discretization.
\\
Finally in \cref{sec:existence_control_problem} the existence of the controls in the time continuous and time discrete case is shown.
In addition the convergence of a subsequence of time discrete  optimal controls to an optimal control of the original problem is obtained.
These results hold not only for aiming at an end time state but also for steering to a state over the whole time horizon.


\section{Solution \revtwo{to} the time-discretized and \revtwo{to} the time-continuous state equations}
\label{sec:analysis_continuous}

First we introduce the time discretization.
Then a certain boundedness property like in \cite{Elliott1996} \revone{(see \eqref{estimates_for_y_tau})}
is shown
which is essential not only for the existence of the solution of the state equation
but also
for proving the existence of an optimal control
and the convergence of the solution of the discretized problem to the time continuous solution.
To obtain this result,
the potential $\psi $
is approximated (as, e.g., in
\cite{Elliott1996}  and \cite{Casas1995})
with a sequence of functions $f_n$ with bounded second derivatives, such that the dominated convergence theorem can be used.
Following the lines of \cite{Elliott1996} we have no restriction on the space dimension $d$.

The existence \revone{for} the time continuous problem will then be shown by taking the limit with respect of the time resolution which also shows convergence of the discretization method.
(In \cite{Elliott1996} first the limit in the time discretization and then in the approximation of $\psi$ is taken).

From now on if no subscripts are provided, with $(\cdot,\cdot)$ and $\|\cdot\|$ we mean the $L^2$- or $\ell_2$-scalar product and norm respectively. The space should be clear from the context. For a Banach space $V$ we will denote its dual by $V'$ and the duality product by $\braket{\cdot, \cdot}$.

\kommentar{{geht dies auch durch mit strictly monotone and coercive?}
 Dacorogna Theorem 3.30 braucht nur, dass $A$ strictly convex und einen bound $A(p) \geq |p|^l$ mit $l >1$. Wir brauchen aber sowieso so etwas wie $c(-1+|p|^2) \leq A(p)$. Zusammen mit der Konvexität ist das wahrscheinlich	 das gleiche wie strongly convex (oder nur marginal verschieden).}

Next we introduce a time discretization and show the existence of a solution of the discretized state equation.
We divide the interval $[0,T]$ into subintervals $I_j \coloneqq (t_{j-1}, t_j]$ for $j = 1,\ldots, N$ with $0=t_0<t_1<\ldots<t_N=T$ and define $\tau_j:=t_j-t_{j-1}$
and $\tau:=\max_j \revone{\tau}_j$.
The state equation we discretize in time with a discontinuous-Galerkin method (dG(0)).
Therefore, let us define
\begin{equation}
\begin{aligned}
\label{eq:def_Y_U}
&  Y_{\tau}:=\{ y_\tau :Q\rightarrow\mathbb{R} \mid
  y_\tau(t,.)\in H^1(\Omega) \revone{\;\forall t}, y_\tau(.,x) \text { \revone{a.e.} constant in } I_j   \text { for } j=1,\ldots,N \} ,\\
 & U_{\tau}:=\{ u_\tau :Q\rightarrow\mathbb{R} \mid
  u_\tau(t,.)\in L^2(\Omega) \revone{\;\forall t}, u_\tau(.,x) \text { \revone{a.e.} constant in } I_j  \text { for } j=1,\ldots,N \} ,
\end{aligned}
\end{equation}
and 
for each interval we label the constant by a subscript, e.g., $y_j:=y_\tau \vert_{I_j}$.
The vector containing these constants will be denoted {by} $(y_j)_{j = 1,\ldots, N} \in H^1(\Omega)^N$.
The time-discretized variant of
\cref{eq1a}--\cref{weakAC}
is then given by
\begin{equation}
\label{discrete_problem1}
\min_{Y_\tau\times U_\tau} {J}({y_\tau},{u_\tau})
= \dfrac{1}{2}\|y_N-y_\Omega\|^2 + \dfrac{\lambda}{2}
\sum_{j=1}^N \tau_j \|u_{j}\|^2
\end{equation}
subject to the time-discretized state equation
\begin{equation}
\label{scheme_state1}
(y_{j}, \varphi_{}) + {\tau_{j}}  (A'(\nabla y_{j}), \nabla \varphi_{}) + {\tau_{j}}(\psi'(y_{j}), \varphi_{}) ={\tau_{j}}(u_{j}, \varphi_{}) + (y_{j-1},\varphi_{}) \quad \forall \varphi \in H^1(\Omega)
\end{equation}
with $j = 1,\ldots,N$ and $y_\tau(0,.):=y_0\in H^1(\Omega)$ is given.

\kommentar{
Here and in the following we use the same notation
for the cost functional, i.e.
we use ${J}$ whether applied to
$(\boldsymbol{y},\boldsymbol{u})\in {H^1(\Omega)^N \times L^2(\Omega)^N} $ or
to $(y_\tau, u_\tau) \in {Y_\tau\times U_\tau}  $ and define
${J}(\boldsymbol{y},\boldsymbol{u}):= {J}(y_\tau, u_\tau)$.
}
We note that the state equation could have arised equally well from an implicit Euler discretization
and we will use the notation $\partial_t^{-\tau} y_\tau$ with 
\begin{equation*}
\partial_t^{-\tau} y_{\tau\vert I_j} \coloneqq \tfrac{1}{\tau_j} (y_j - y_{j-1})
\end{equation*}
\revtwo{in $L^2(\Omega)$ for $j=1,\ldots,N$.}\\
One may favour a splitting approach for $\psi$ or an approximation of the quasilinear term $A$ as in \cite{bgn13, Barrett2014} instead of the fully implicit method. However, to our knowledge there exists no convergence proof for these discretizations of the state equation {to the time continuous one in the limit $\tau \to 0$}.
{Moreover, while for \johannes{the} implicit time discretization \johannes{the} differentiability  of the control to state operator
is obtained in   \cite{BlMeNumeric}
under additional smoothness properties on
$A$ and $\psi$,
  it is \revone{not known whether} this property holds for semi-implicit discretizations.
The additional computational cost using implicit discretization} is nearly negligible for solving the optimal control problem.\\
The first step is given by the subsequent lemma.
\begin{lemma}\label{lemma1}
Let $A$ fulfill the conditions in \cref{assu}.
Furthermore let $y_0 \in H^1(\Omega)$ and
$u_\tau \in U_\tau$.
\revtwo{Let $f\in C^2(\mathbb{R})$ be a function such that $\vert f'' \vert\in L^\infty(\mathbb{R})$ and $f''\geq -C_\psi$ on $\mathbb{R}$ for some constant $C_\psi\geq0$. Then, for any $0<\tau <\tfrac{1}{C_\psi}$, there exists a function $y_\tau  \in Y_\tau$ which is} a solution of
$	y_\tau(0)=y_0 $ in  $\Omega$ and for all $j=1,\ldots, N$ it holds
$\forall \eta \in  H^1(\Omega): $
\begin{equation}
\begin{gathered}
\label{proof_existence_bounded}
\int_\Omega \partial^{-\tau}_t \revone{y_{\tau\vert I_j}} (\eta-y_j) + A(\nabla \eta)-  A(\nabla y_j) + f'(y_j)(\eta-y_j) -u_j(\eta-y_j) \geq 0.
\end{gathered}
\end{equation}
\revtwo{In addition, if $\Lambda>0$ and $u_\tau, y_0, A, f$ fulfill that}
\\
$\norm{u_\tau}_{L^2(0,T; L^2(\Omega))}, 
\norm{y_0}_{L^2(\Omega)} \leq \Lambda$,
        $-\Lambda+\Lambda^{-1}\vert p\vert^2 \leq A(p) \text{ and } A'(p)^Tp \leq \Lambda \vert p \vert^2$ as well as
\begin{equation}
\label{f_lambda_bounds}
\int_{\Omega}(A(\nabla y_0) + f(y_0)) \leq \Lambda \quad \text{ and } \quad
f\geq -\Lambda  ,  \quad f''\geq -\Lambda,
\end{equation}
\revtwo{then there exist constants $\tau_\Lambda>0$ and $C(\Lambda)>0$, depending on $\Lambda$, and for all $0<\tau<\tau_\Lambda$, the solutions $y_\tau$ of \eqref{proof_existence_bounded} satisfy that}
\begin{equation}
\label{estimates_for_y_tau}
\begin{gathered}
\norm{\partial_t^{-\tau}y_\tau}_{L^2(0,T;L^2(\Omega))}+ \norm{y_\tau}_{L^\infty(0, T; H^1(\Omega))}
+ \norm{f'(y_\tau)}_{L^2(0, T; L^2(\Omega))}
\leq C(\Lambda).
\end{gathered}
\end{equation}
\end{lemma}
Note that under \cref{assu} one can find for given  $u_\tau,y_0,A$ and $f$ always a constant $\Lambda$ such that the above required estimates hold. In particular for $A$
the growth condition induces \revone{$A'(p)^Tp \leq \overline{C}_A \vert p \vert^2$ as well as $A'(0)=0$ and then the strong monotonicity provides
$A(0)+\tfrac{1}{2}C_A\vert p\vert^2 \leq A(p)$.}
\begin{proof}
  We note that $f$ and $f'$ induce \revone{continuous} Nemytskii operators $f:L^2(\Omega) \to L^1(\Omega)$ and
  $ f':L^2(\Omega) \to L^2(\Omega)$ due to the bounds on $f''$.\\
Starting with $y_{0}:=y(0)$,
 define $y_j\in H^1(\Omega)$  successively for $j\geq 1$
	to be the unique minimizer of
 \begin{equation}\label{minstat}
\begin{gathered}
  \Phi_{j,\tau}(\eta) \coloneqq
  \int_{\Omega} \left( \tfrac{1}{2\tau_j} \vert \eta - y_{j-1} \vert^2 + A(\nabla \eta) + f(\eta) - u_j\eta \right)
	\end{gathered}
\end{equation}
where the integrands are strongly convex for $\tfrac1\tau +f''(s) \geq \tfrac1\tau -C>0$.
\revone{For $\eta \in H^1(\Omega)$, $\delta > 0$
we obtain with $\eta_\delta \coloneqq y_j + \delta(\eta - y_j)$ and using
the convexity of $A$ as well as
$\Phi_{j,\tau}(y_j) \leq \Phi_{j,\tau}(\eta_\delta)$
\begin{equation}
\label{eq:tmp_VI_for_disc}
\begin{aligned}
&  \int_\Omega (A(\nabla \eta) - A(\nabla y_j) )
\geq \int_\Omega \dfrac{1}{\delta} (A(\nabla \eta_\delta) - A(\nabla y_j))\\
&~~	\geq - \int_\Omega (\tfrac{1}{2\tau_j\delta} (\vert\eta_\delta - y_{j-1}\vert^2 - \vert y_j-y_{j-1}\vert^2) + \tfrac{1}{\delta}(f(\eta_\delta) - f(y_j)) - u_j(\eta-y_j) ) \\
&~~	\stackrel{\delta \to 0}{\to} -\int_\Omega \left(\dfrac{y_{j} - y_{j-1}}{\tau_j}\right)(\eta - y_j) + f'(y_j)(\eta - y_j) - u_j(\eta-y_j).
\end{aligned}
\end{equation}
The last term is obtained using the mean value theorem and applying dominated convergence given $\vert f'(x)\vert \leq C(1+\vert x\vert)$ (see also \cref{rem:take_limit_in_VI}). Altogether we have shown \cref{proof_existence_bounded}.}
Now we want to deduce the estimate \eqref{estimates_for_y_tau}.
The summation of $\Phi_{l,\tau}(y_l) \leq \Phi_{l,\tau}(y_{l-1})$ 
yields
\kommentar{
	\begin{align*}
	&\sum_{i=1}^{j} \int_{\Omega} \dfrac{1}{2\tau_i}\left|y_i - y_{i-1}\right|^2 + A(\nabla y_i) + f(y_i) &\leq& \sum_{i=1}^{j} \int_{\Omega} A(\nabla y_{i-1}) + f(y_{i-1}) + u_i(y_i-y_{i-1}) \\&&=& \sum_{i=0}^{j-1} \int_{\Omega} A(\nabla y_{i}) + f(y_{i}) + \sum_{i=1}^{j} \int_{\Omega} u_i(y_i-y_{i-1})
	\end{align*}
	\begin{equation*}
	\iff \quad \underbrace{\sum_{i=1}^{j} \int_{\Omega} \dfrac{1}{2\tau_i}\left|y_i - y_{i-1}\right|^2}_{=\tfrac{1}{2}\int_{0}^{t_j}\int_\Omega |\partial_t^{-\tau} y_\tau|^2} + \int_\Omega A(\nabla y_j) + f(y_j) \leq \int_{\Omega} A(\nabla y_{0}) + f(y_{0}) + \underbrace{\sum_{i=1}^{j} \int_{\Omega} u_i(y_i-y_{i-1}}_{=\int_{0}^{t_j}\int_{\Omega} u_\tau \partial_t^{-\tau} y_\tau}),
	\end{equation*}
	and finally using (scaled) Young's inequality
      }%
\begin{align}
&	\int_{0}^{t_j} \int_{\Omega} \tfrac{1}{2} \vert \partial_t^{-\tau} y_\tau\vert^2 + \int_{\Omega} (A(\nabla y_j) + f(y_j)) \nonumber \\
&~~\leq \int_{\Omega} (A(\nabla y_0) + f(y_0)) + 	\int_{0}^{t_j} \int_{\Omega} u_\tau \partial_t^{-\tau} y_\tau
\leq C(\Lambda) + \tfrac{1}{4} \int_{0}^{t_j} \int_{\Omega} \vert \partial_t^{-\tau} y_\tau\vert^2.
\label{eq:9}
\end{align}
Using the assumptions $-\Lambda+\Lambda^{-1}\vert p\vert^2 \leq A(p)$, $- f\leq \Lambda$,
$\norm{y_0}_{L^2(\Omega)}\leq \Lambda $ as well as $y_j= y_0 + \int_{0}^{t_j} \partial_t^{-\tau}y_\tau$,
we obtain
\kommentar{
  see https://math.stackexchange.com/questions/1104652/bounded-second-derivative-for-a-convex-function-implies-growth-bound-on-the-funt
  }
  \begin{equation}
    \label{eq:estimates_y_tau_discrete_in_continuous_formulation}
	\begin{gathered}
	\norm{\partial_t^{-\tau}y_\tau}_{L^2(0,T;L^2(\Omega))}+  \norm{\nabla y_\tau}_{L^\infty(0,T; L^2(\Omega))}
	+ \norm{y_\tau}_{L^\infty(0, T; L^2(\Omega))} \leq C(\Lambda).
	\end{gathered}
      \end{equation}
Then, choosing $\eta\coloneqq y_{j} - \delta f'(y_{j}), \delta >0$ in \cref{eq:tmp_VI_for_disc}, we obtain
\begin{equation}
\label{proof_Elliot_1}
\int_{\Omega} f'(y_{j})^2
\leq \int_{\Omega} - \tfrac{y_{j} - y_{j-1}}{\tau_j}f'(y_{j}) + u_j f'(y_{j}) - \dfrac{1}{\delta} (A(\nabla y_{j}) - A(\underbrace{\nabla(y_{j}-\delta f'(y_{j}))}_{= \revone{\nabla y_j-\delta f''(y_j)\nabla y_j}})).
\end{equation}
To the third integral we \revone{applied} the mean value theorem pointwisely almost everywhere \revtwo{in $\Omega$}, with the intermediate point of 1 and $1-\delta f''(y)$ denoted by $\xi_\delta(y)$.
Note that due to the boundedness of $f''$ also $\xi_\delta(\cdot)$ is bounded and $\xi_\delta \to 1$ \revtwo{as} $\delta \to 0$, \revtwo{in the pointwise sense}. Now we can use $0 \leq A'(p) ^T p \leq \Lambda\vert p\vert^2$
and $-f''\leq \Lambda$ as well as dominated convergence to obtain
\begin{eqnarray*}
		\int_{\Omega} f'(y_{j})^2
		&\leq& \int_{\Omega} - \tfrac{y_{j} - y_{j-1}}{\tau_j}f'(y_{j}) + u_j f'(y_{j}) + C(\Lambda) \vert\nabla y_{j}\vert^2\\
		&\leq& \int_{\Omega} \left(\tfrac{y_{j} - y_{j-1}}{\tau_j}\right)^2 + \tfrac{1}{4}f'(y_{j})^2 + \tfrac{1}{2} u_j^2 + \tfrac{1}{2}f'(y_{j})^2 + C(\Lambda) \vert\nabla y_{j}\vert^2,
	\end{eqnarray*}
and hence  with \cref{eq:estimates_y_tau_discrete_in_continuous_formulation}
\begin{equation}
\label{eq:estimates_for_y_disc}
\left( \sum_{j=1}^N \tau_j \norm{f'(y_{j} )}^2 \right) ^{1/2}\leq C(\Lambda).
\end{equation}
\end{proof}
\begin{remark}
  \label{rem:take_limit_in_VI}
 As \revtwo{was} mentioned, the strong monotonicity of $A'$ and the growth condition
  $ \vert A'(p)\vert \leq \overline{C}_A\vert p\vert $ induces
  $A'(0)=0$,
 $A'(p)^Tp \leq \overline{C}_A \vert p\vert^2$ and
 $A(0)+\tfrac{1}{2}C_A\vert p\vert^2 \leq A(p)  $. Furthermore
 $A(p) \leq c(1+\vert p\vert^2)$ \revtwo{also} holds with some $c>0$. 
 Hence \revtwo{for all  $\eta,\xi \in  L^2(0,T;H^1(\Omega))$ we have}
  $A(\nabla \eta)\in L^1(Q)$ and
  (using Young's inequality) $A'(\nabla \eta)^T\nabla \xi \in L^1(Q)$.
 It also induces the pointwise estimate
  $ \vert A'(\nabla y + s \delta \nabla \xi)^T\nabla \xi\vert \leq C(\vert\nabla y\vert^2 + \vert\nabla \xi\vert^2),
  $ for $0\leq s \delta \leq 1$
providing  an integrable majorant, which allows
to take the limit 
$\delta \searrow 0$ for the integral below.
Hence we obtain
\begin{equation} \label{ineqtoeq}
\begin{aligned}
  \lim_{\delta \searrow 0}\dfrac{1}{\delta} \int_Q (A(\nabla y + \delta \nabla \xi) - A(\nabla y) )
  &=
      \lim_{\delta \searrow 0}\int_Q \int_{0}^{1} A'(\nabla y+s \delta \nabla \xi)^T \nabla \xi \, \mathrm ds
      \\
      &= \int_Q A'(\nabla y)^T \nabla \xi.
\end{aligned}
\end{equation}
The same holds respectively for integration over $\Omega$.
 Together with the monotonicity of $A'$ this enables the usual steps of the proof that solving
 the variational inequality
\cref{proof_existence_bounded}
is equivalent to solving the variational equality
\cref{scheme_state1}
with $f$ instead of $\psi$. 
  \end{remark}

With \cref{lemma1} at hand
we can show the existence of a unique weak solution to the time discretized state equation \cref{scheme_state1}.
{Note that the following bound \cref{estim_for_y_tau} (and likewise \cref{eq:bounds_on_y_in_lemma_limit_from_discrete} in the time-continuous case) will be crucial for showing the existence of an optimal control later.}
\begin{theorem}
\label{th:existence_state_disc}
Let \cref{assu} be fulfilled.
If $\tau =\max_j \tau_j < 1/C_\psi$
then for every $u_\tau \in U_\tau$
the time discretized state equation \cref{scheme_state1} has a unique solution $y_\tau \in Y_\tau$.\\
The solution operator is denoted by $S_\tau:U_\tau\rightarrow  Y_\tau$.
\\
Furthermore, \revtwo{if $\norm{u_\tau}_{L^2(0,T;L^2(\Omega))}\leq \bar c$ for all $0<\tau<1/C_\psi$, then there exists a constant $C_{A,\psi,y_0}(\bar{c})$, independent of $\tau$, such that}
\begin{equation}
\label{estim_for_y_tau}
\begin{gathered}
\norm{\partial_t^{-\tau}y_\tau}_{L^2(0,T;L^2(\Omega))}+ \norm{y_\tau}_{L^\infty(0, T; H^1(\Omega))}
+ \norm{\psi'(y_\tau)}_{L^2(0, T; L^2(\Omega))}
\leq C_{A,\psi,y_0}(\bar c).
\end{gathered}
\end{equation}
\end{theorem}
\begin{proof}
We consider the approximation of $\psi$ by $f_n$ according to Assumption \ref{assu}. Then,  
 due to $\psi(y_0)\in L^1(\Omega)$ for given $y_0\in H^1(\Omega)$ and $-c\leq f_n\leq c(\psi+1)$, $ -C_\psi\leq f_n''$ 	one can find $\Lambda$\revone{,} depending only on ${A,\psi,y_0}$\revone{,} large enough such that
\begin{equation}
\begin{gathered}
\label{condlemma}
\norm{y_0}_{L^2(\Omega)}, \, \int_{\Omega}(A(\nabla y_0) +{f_n}(y_0)), \,  -\inf_{t\in \mathbb{R}} f_n(t), \, -\inf_{t\in \mathbb{R}} f_n''(t) \leq \Lambda,\\
	-\Lambda+\Lambda^{-1}\vert p\vert^2 \leq A(p) \text{ and } A'(p)^Tp \leq \Lambda \vert p\vert^2.
	\end{gathered}
      \end{equation}
 We denote by $y _{j,n}$ the solutions	 of
 \cref{proof_existence_bounded} with $f=f_n$
 which exist according to \cref{lemma1} and remark that they exist  for $ \tau<\tfrac1{C_\psi}$
 where the integrands of
 \cref{minstat}
 are strongly convex due to $\tfrac1\tau +f_n''(s) \geq \tfrac1\tau -C_\psi>0$.
 Also \cref{lemma1} provides the estimates
 \cref{estimates_for_y_tau}, i.e., for all $\tau $ and $n$ it holds
 \begin{equation}
\label{estimay_taun}
\begin{gathered}
\norm{\partial_t^{-\tau}y_{\tau,n}}_{L^2(0,T;L^2(\Omega))}+ \norm{y_{\tau,n}}_{L^\infty(0, T; H^1(\Omega))}
+ \norm{f_n'(y_{\tau,n})}_{L^2(0, T; L^2(\Omega))}
\leq C(\Lambda,\bar c).
\end{gathered}
\end{equation}
Then\revtwo{, \cite[Lemma 1.3]{Lions:233038} together with} $f_n \to \psi$ in $C^1_{\text{loc}}$ for $n \to \infty$ and
the weak-lower semicontinuity of \revtwo{$A: w\in L^2(\Omega;\mathbb{R}^d)\mapsto \int_\Omega A(w)\, \mathrm dx \in \mathbb{R}$}
allows to take for a subsequence of $y_{j,n} $ the \revtwo{limit inferior for} $n\to\infty$ for all terms in \cref{estimay_taun} and
\revtwo{the limit superior for $n\to\infty$ for those in}
\cref{proof_existence_bounded} 
to obtain {\cref{estim_for_y_tau}} and
that for all $\eta \in H^1(\Omega )$ it holds
\begin{equation}
\label{eq:ineqwithpsi}
\int_\Omega (A(\nabla \eta) - A(\nabla y_{j}) )
	\geq  \int_\Omega - \left(\dfrac{y_{j} - y_{j-1}}{\tau_j}\right)(\eta - y_j) - \psi'(y_j)(\eta - y_j) + u_j(\eta-y_j) .
	\end{equation}
	Finally we can go over to the equality \cref{scheme_state1}
	by the reasoning from \cref{rem:take_limit_in_VI}.
        
	The uniqueness of the solution of  \cref{scheme_state1} can be shown for each time step separately one after another. For this purpose assume the existence of two solutions.
	Subtracting their defining equations, testing with their difference and using the strong monotonicity of $A'$ and
	of $s+\tau_j\psi'(s) $
	due to $\tau < 1/C_\psi$
	shows that the $H^1$-norm of their difference vanishes.
\end{proof}
With a further (minor) restriction on the maximal time step $\tau$ we obtain Lipschitz-continuity \revtwo{of the solution operator for \eqref{scheme_state1} with a constant independent \revone{of} $\tau$}.
\begin{theorem}
	\label{th:lipschitz_disc}
\begin{samepage}
		Let \cref{assu} 
	and
	$\tau \leq \tfrac{1}{1+2C_\psi}$ hold.
	Then the mapping $\tilde S_\tau:(y_0, u_\tau) \mapsto y_\tau$
	where  $ y_\tau$ is the solution of equation \cref{scheme_state1}, is Lipschitz-continuous \revtwo{in the sense that} 
	\begin{equation}
		\begin{aligned}
		\label{eq:Lipschitz_discrete_inequ2}
		\norm{y^{(1)}_\tau-y^{(2)}_\tau&}_{L^\infty(0,T;L^2(\Omega))} + \norm{\nabla y^{(1)}_\tau-\nabla y^{(2)}_\tau}_{L^2(0,T;L^2(\Omega))} \leq\\
		&\leq {C}_{A,\psi,T}\left(\norm{y^{(1)}_0 - y^{(2)}_0}_{L^2(\Omega)} + \norm{u^{(1)}_\tau - u^{(2)}_\tau}_{L^2(0,T;H^1(\Omega)')}\right),
		\end{aligned}
	\end{equation}
	where $y_{\tau}^{(i)}= \tilde S_\tau(y_0^{(i)}, u_\tau^{(i)}) $ for $i=1,2$.
\end{samepage}
\end{theorem}
\begin{proof}
	We note down the differences by a prescript $\delta$, e.g., $\delta y_\tau \coloneqq y_{\tau}^{(1)} - y_{\tau}^{(2)}$.
	With $\tfrac{1}{2}(a^2-b^2) \leq (a-b)a$ \revtwo{in mind}, 
	testing the defining equalities \cref{scheme_state1} with $\delta y_j$ and using that $A'$ is strongly monotone as well as \cref{eq:lower_estimate_psi_prime_diff}, we obtain
	\begin{align*}
	\tfrac{1}{2}&\left(\norm{\delta y_j}^2 - \norm{\delta y_{j-1}}^2\right) + \tau_j C_A \norm{\nabla \delta y_j}^2\\ &\leq \left(\delta y_j - \delta y_{j-1}, \delta y_j\right) + \tau_j\left(A'(\nabla y_{j}^{(1)})-A'(\nabla y_j^{(2)}), \nabla \delta y_j \right) \\
	&= \tau_j \left(\delta u_j, \delta y_j\right) - \tau_j \left(\psi'(y_{j}^{(1)}) -\psi'(y_j^{(2)}), \delta y_j\right)\\
	&\leq \tfrac{\tau_j}{2\epsilon} \norm{\delta u_j}_{{H^1}'}^2 + \tfrac{\tau_j\epsilon}{2}\norm{\delta y_j}_{H^1}^2 + \tau_jC_\psi\norm{\delta y_j}^2.
	\end{align*}
	In the last step we used scaled Young's inequality with $0<\epsilon <\min(1, 2C_A)$.%
	\kommentar{$\epsilon<1$ braucht man für $C_{\psi, \tau}>0$ später, damit nur $\tau \leq \ldots$ und nicht $\tau<\ldots$ gefordert ist. Das wiederum wird benötigt, um ohne Probleme $C_{\psi,\tau}$ unabhängig von $\tau $ abzuschätzen.}
	We now sum over $j = 1,\ldots,J$ and get
	\begin{align}
	\label{eq:discr_cont_Ausgangsgleichung}
	\tfrac{1}{2}\norm{\delta y_J}^2 + \tilde{C}_A \sum_{j=1}^{J} {\tau_j} \norm{\nabla \delta y_j}^2 \leq \tfrac{1}{2} \left(\norm{\delta y_0}^2 + \sum_{j=1}^{J} \tfrac{\tau_j}{\epsilon} \norm{\delta u_j}_{{H^1}'}^2\right) + \tfrac{1}{2}\tilde{C}_\psi \sum_{j=1}^{J} \tau_j \norm{\delta y_j}^2
	\end{align}
	for all $1 \leq J \leq \revtwo{N}$. Here we defined $\tilde{C}_A\coloneqq C_A-\tfrac{\epsilon}{2}$ and $\tilde{C}_\psi \coloneqq {\epsilon}+2C_\psi$
	Omitting the gradient term on the left and absorbing the $J$-th term from the right, we obtain
	\begin{align*}
	\norm{\delta y_J}^2 &\leq \dfrac{1}{
		(1-\tilde{C}_\psi\tau_J)} \left(\norm{\delta y_0}^2 + \sum_{j=1}^{J} \tfrac{\tau_j}{\epsilon} \norm{\delta u_j}_{{H^1}'}^2\right) + \dfrac{\tilde{C}_\psi}{1-\tilde{C}_\psi\tau_J}\sum_{j=1}^{J-1} \tau_j \norm{\delta y_j}^2\\
	& \leq C_{\psi,\tau} \left(\norm{\delta y_0}^2 + \sum_{j=1}^{\revtwo{N}} \tfrac{\tau_j}{\epsilon} \norm{\delta u_j}_{{H^1}'}^2\right) +
	C_{\psi,\tau}{\tilde{C}_\psi}\sum_{j=1}^{J-1} \tau_j \norm{\delta y_j}^2,
	\end{align*}
	where \revtwo{we had to suppose smallness of $\tau$ to require}
	$C_{\psi,\tau}:= \tfrac{1}{1-\tilde{C}_\psi\tau} >0$.
	\revtwo{Now} we apply the discrete Gronwall Lemma
	\revtwo{(see, e.g. \cite[Lemma A.3]{Kruse2014}),}
	which yields
	\begin{equation}
	\label{eq:th_disc_Gronwall_est_1}
	\begin{aligned}
	\norm{\delta y_J}^2 &\leq \left(\norm{\delta y_0}^2 + \sum_{j=1}^{\revtwo{N}} \tfrac{\tau_j}{\epsilon} \norm{\delta u_j}_{{H^1}'}^2\right) C_{\psi,\tau} \exp\left(C_{\psi,\tau}{\tilde{C}_\psi} \sum_{j=1}^{J-1} \tau_j\right)\\
	&\leq \left(\norm{\delta y_0}^2 + \sum_{j=1}^{\revtwo{N}} \tfrac{\tau_j}{\epsilon} \norm{\delta u_j}_{{H^1}'}^2\right)
	C_{\psi,\tau} \exp\left(C_{\psi,\tau}{\tilde{C}_\psi T}\right).
	\end{aligned}
	\end{equation}
	Inserting this into \cref{eq:discr_cont_Ausgangsgleichung} we finally get for all $J=1,\ldots,\revtwo{N}$
	\begin{equation}
	\label{eq:th_disc_Gronwall_est_2}
	\tilde{C}_A \sum_{j=1}^{J} {\tau_j} \norm{\nabla \delta y_j}^2 \leq \tfrac{1}{2} \left(\norm{\delta y_0}^2 + \sum_{j=1}^{\revtwo{N}} \tfrac{\tau_j}{\epsilon} \norm{\delta u_j}_{{H^1}'}^2\right)
	\left(1 + C_{\psi,\tau}{\tilde{C}_\psi T}
	\exp\left(C_{\psi,\tau}{\tilde{C}_\psi T}\right)\right),
	\end{equation}
        which together with \cref{eq:th_disc_Gronwall_est_1}
and the boundedness of $C_{\psi,\tau}$ independently of $\tau$
        yields the inequality \cref{eq:Lipschitz_discrete_inequ2}. 
\end{proof}
A similar result \revtwo{to} \cref{th:existence_state_disc} could also be obtained by using results on monotone operators, see, e.g., \cite{Lions:233038}. 
Together with an argument formerly found by Stampacchia one would obtain the regularity $y_j\in L^\infty(\Omega) \cap H^1(\Omega)$ at each step of our time discretization \cite{Stampacchia1965, Casas1995}.
These results are applicable if $\tau$ is sufficiently small such that the term $y_j + \tau_j \psi'(y_j)$ becomes monontonic. However this regularity comes with 
restriction on the space dimension $d$.
\kommentar{Brauchen wir das?: In contrast to this we are able to show that the norm $\norm{\psi'(y_j)}_{L^2(\Omega)}$ is bounded independently from the space dimension. When considering the first order conditions for the discretized problem in our follow-up paper, we will however have to use the $L^\infty$-bound and restrict $d$.}

Our approach also allows taking the limit $\tau\to 0$ providing the \johannes{convergence of the time discrete solutions $y_\tau$ with $\tau\to0$ and with it} the existence of a solution in the time-continuous case. 
The following lemma serves as a preparation and gives a similar result to \cite{Dreher2012}. 

\begin{lemma}\label{th:conv_disc_to_cont_solution}
	Let $\{y_\tau\}_{\tau}$ with $y_\tau\in Y_\tau$ and $\tau\to0$ be a sequence satisfying
	\begin{equation}
		\label{eq:ass_for_C_conv}
		\norm{\partial_t^{-\tau}y_\tau}_{L^2(0,T;L^2(\Omega))}+ \norm{y_\tau}_{L^\infty(0, T; H^1(\Omega))} \leq C,
	\end{equation}
	where $C>0$ is independent of $\tau$. 
	Then there exists a subsequence (again denoted by $\{y_\tau\}_{\tau}$) and a function $z\in C([0,T]; L^2(\Omega))$, such that
	\begin{equation}
		\label{eq:y_t_strong_convergence_final}
		\max_{t\in[0,T]} \| y_\tau(t) -z(t)\|_{L^2(\Omega)} \to 0 \qquad \text{as} \qquad \tau\to0.
	\end{equation}
\end{lemma}
\begin{proof}
	Using the definitions from \cref{eq:def_Y_U}, for given $y_\tau$ we define its linear interpolant $z_\tau$, 
	i.e.,
	$		z_\tau(t)_{\vert I_{j}} = y_{j-1} + (t-t_{j-1})\partial_t^{-\tau} y_\tau(t_j).
	$
	Note that from \cref{estim_for_y_tau} we have $\|z_\tau\|_{H^1(0,T;L^2(\Omega))\cap L^\infty(0,T;H^1(\Omega))}\leq C$ \revtwo{for some constant $C$,} independent \revtwo{of}~$\tau$. By the compact imbedding $L^\infty(0,T;H^1(\Omega))\cap H^1(0,T;L^2(\Omega)) \hookrightarrow C([0,T];L^2(\Omega))$
	(see Aubin--Lions--Simon compactness  theorem%
	, e.g.,
	in
	\cite{Simon1986})
	we deduce the existence of a \revtwo{function} $z$ such that (possibly for a subsequence) $z_\tau\to z$ in $C([0,T]; L^2(\Omega))$.
	In addition, for 
	$t = \beta t_j + (1-\beta) t_{j-1}$ with $\beta\in(0,1]$ we find
	\begin{equation}
		\begin{aligned}
			\|y_\tau(t) - z_\tau(t)\|_{L^2(\Omega)}^2 
			&
			=(1-\beta)^2(t_j-t_{j-1})^2\| \partial_t^{-\tau} y_\tau(t_j) \|_{L^2(\Omega)}^2\\
			&\leq (1-\beta)^2 \tau \|\partial_t^{-\tau} y_\tau\|_{L^2(0,T;L^2(\Omega))}^2
			\leq C_{A,\psi,y_0} \tau
		\end{aligned}
	\end{equation}
	independent of $t$.
	Consequently it holds 
	$\max_{t\in[0,T]} \| y_\tau(t) -z(t)\|_{L^2(\Omega)} \to 0 
	$ as $\tau \to 0$.
\end{proof}

We are now prepared to deduce the \johannes{convergence of the numerical solutions to} the time-continuous \johannes{solution}.
\begin{theorem}
	\label{existence_state}~\\
Let \cref{assu} hold.
Then for every $u \in L^2(0,T;L^2(\Omega))$ there exists a unique weak solution
$y\in L^\infty(0,T;H^1(\Omega)) \cap H^1(0,T;L^2(\Omega))$
to \cref{weakAC},
i.e.,
\begin{equation*}
\int_Q  \partial_t y \eta + A'(\nabla y)^T \nabla \eta +  \psi'(y) \eta =  \int_Q u\eta 
\qquad  \forall \eta \in L^2(0, T; H^1(\Omega)),
\end{equation*}
\revtwo{subject to $y(0)=y_0\in H^1(\Omega)$ a.e. in $\Omega$, and}
{it also holds $\psi'(y)\in L^2(Q)$.}\\
For a sequence $\{u_\tau\}$ with $u_\tau \rightharpoonup  u$ in $L^2(0,T;L^2(\Omega))$ for $\tau\to 0$, then for the corresponding sequence of time-discrete solutions $\{y_\tau\}$ and $y\coloneqq S(u)$ the following convergences hold:
\begin{equation}
\label{eq:various_convergences_of_y_tau}
\begin{aligned}
y_\tau\rightharpoonup y \text{ in } L^2(0,T;H^1(\Omega)),\\y_\tau\stackrel{\ast}{\rightharpoonup} y \text{ in } L^\infty(0,T;H^1(\Omega)),\\
 y_\tau \to y \text{ in } C([0,T];L^2(\Omega)),\\ \partial^{-\tau}_t y_\tau \rightharpoonup \partial_t y \text{ in } L^2(0,T;L^2(\Omega)),\\ \psi'(y_\tau) \rightharpoonup \psi'(y)  \text{ in } L^2(0,T;L^2(\Omega)).
\end{aligned}
\end{equation}
\end{theorem}

\begin{proof}
Given $u\in L^2(0,T;L^2(\Omega))$ we choose a sequence of discretizations $u_\tau \in U_\tau$ with $u_\tau \rightharpoonup  u$ in $L^2(0,T;L^2(\Omega))$ for $\tau\to 0$.
This allows for the choice of a \revtwo{constant} $\bar c >0$ \revtwo{satisfying} $\norm{u_\tau} \leq \bar c$.
Let $y_\tau$ be the solution of \cref{scheme_state1} corresponding to $u_\tau$. Then the estimates
\cref{estim_for_y_tau} hold.
Hence, for $\tau\to 0$ there exists a (sub-)sequence satisfying the first and last two convergences in \cref{eq:various_convergences_of_y_tau}.
\kommentar{es folgt $y_\tau \stackrel{*}{\rightharpoonup} y$ (schwacher Grenzwert gleich $y$ da punktweiser Grenzwert aus anderen Konvergenzen in \cref{eq:various_convergences_of_y_tau} bekannt). Die weak-*-lowersemicontinuity der Norm liefert das Ergebnis.}%
\revone{The strong convergence of $\{y_\tau\}$ in $C([0,T];L^2(\Omega))$ is obtained from 
\cref{th:conv_disc_to_cont_solution}.
Note that for later convenience we already denote the limit as $y$ and the latter two limits in \cref{eq:various_convergences_of_y_tau} were identified using pointwise almost-everywhere convergence of $y_\tau$ (following from the strong convergence), continuity of $\psi'$ and an application of \cite[Lemma 1.3]{Lions:233038}.
Altogether} we can take the limit
in the variational inequality \cref{eq:ineqwithpsi} to obtain that $y$ satisfies
\begin{equation}
  \label{proof_existence_state1}
\int_Q \partial_t y(\eta-y) + A(\nabla \eta)-  A(\nabla y) + \psi'(y)(\eta-y) -u(\eta-y) \geq 0
\quad \forall  \eta \in L^2(0,T;H^1(\Omega)).
\end{equation}
\Cref{rem:take_limit_in_VI} yields that $y$ solves also the variational equality \cref{weakAC}.
Furthermore, using weak ($\ast$) lower-semicontinuity
the solution $y$  satisfies
\begin{equation}
\label{eq:bounds_on_y_in_lemma_limit_from_discrete}
\norm{\partial_t y}_{L^2(0,T;L^2(\Omega))} + \norm{y}_{L^\infty(0,T;H^1(\Omega))} + \norm{\psi'(y)}_{L^2(0,T;L^2(\Omega))}\leq C_{A,\psi,y_0}(\bar c).
\end{equation}
{The uniqueness follows by subtracting the defining equations for two solutions and using a Gronwall argument to deduce that their difference vanishes.}

Recall that the choice of discretization (given by the choice of the intervals) was arbitrary. Furthermore on the way to obtain $y$ from $y_{j,n}$ we had to take subsequences twice. Whatever choice made we would have got a $y$ satisfying the same variational inequality. \revtwo{Since} this  variational inequality has a unique solution, the whole sequence has to converge.
Summarized, for all discretization, we get a sequence $y_{j,n}$ that for $n\to \infty$ and then $\tau\to 0$ ($j \to \infty$) results in the same limit $y$ satisfiying the  variational inequality.
\end{proof}

Finally, we 
\johannes{also obtain}
Lipschitz-continuity of the time-continuous solutions.
\begin{theorem}
\begin{samepage}
		Let \cref{assu} hold.
	Then the solution of \eqref{weakAC} depends Lipschitz-continuously on $(y_0, u)$ \revtwo{in the sense that}
	\begin{equation}
	\label{eq:Lipschitz_dependency_estimate}
	\begin{aligned}
	&	\norm{y_1-y_2}_{C([0,T];L^2(\Omega))\cap L^2(0,T;H^1(\Omega))}\\
	&\quad 	\leq
	C_{\psi, A, T}  \left(\norm{y_{1,0}-y_{2,0}}_{L^2(\Omega)} +
	\norm{u_1-u_2}_{L^2(0,T;H^1(\Omega)')}\right),
	\end{aligned}
	\end{equation}
	where $y_1,y_2$ are the solutions to the data $(y_{1,0}, u_1)$ and $(y_{2,0}, u_2)$ respectively.
\end{samepage}	
\end{theorem}
\begin{proof}
  \johannes{
Using sequences $u_\tau^{(1)}$, $u_\tau^{(2)} \in U_\tau$ converging to $u_{1}$ and $u_{2}$ in $L^2(0,T;L^2(\Omega))$
we obtain a sequence $y_\tau^{(i)} \coloneqq S_\tau(u_\tau^{(i)})$ with $y_\tau^{(i)}$ converging to
$y_{i}= S(u_{i})$ in the sense of \eqref{eq:various_convergences_of_y_tau} for $i=1,2$.
Then applying \cref{th:lipschitz_disc} one obtains the Lipschitz estimate~\eqref{eq:Lipschitz_dependency_estimate} since the constant in \eqref{eq:Lipschitz_discrete_inequ2} is independent of $\tau$.}
\end{proof}

\johannes{We remark that this result can be obtained in a similar way as for the discrete solutions using the continuous Gronwall inequality.}

\revone{In addition, let us mention that our discretization of the state equation inherits an important property of the time continuous case. Namely, since the state equation is a result of the gradient flow of the energy $\mathcal{E}$ given in \cref{GinzburgLandauEnergyAniso}, the energy should decrease in time when there is no input, i.e., $u=0$.}
\begin{theorem}
	\label{th:energy_stable_disc}
  Let \cref{assu}
  and $\tau \leq {2}/{C_\psi}$
  hold.
  Then the scheme  \cref{scheme_state1} for the state equation is energy stable, i.e., for $u_\tau =0$
  the energy functional $\mathcal{E}$ is decreasing in time.
\end{theorem}
\begin{proof}
We test \cref{scheme_state1} with the difference $y_{j}-y_{j-1}$ and obtain
\begin{equation}
  \dfrac{1}{\tau_j}\norm{y_{j} - y_{j-1}}^2 + \left(A'(\nabla y_j), \nabla y_{j} - \nabla y_{j-1}\right) +
  \left(\psi'(y_{j}), y_{j}- y_{j-1}\right) = 0.
 \end{equation}
The convexity of $A$ (recall $A'$ is strongly monotone) yields
$ \left(A'(\nabla y_j), \nabla y_{j} - \nabla y_{j-1}\right) \geq A(\nabla y_{j}) - A(\nabla y_{j-1})$
for the second term. The third term can be estimated by the following relation
\begin{equation}\label{eq:psitaylor}
\psi'(y_{j})( y_{j}- y_{j-1}) \geq \psi(y_{j}) - \psi(y_{j-1}) - \tfrac{C_\psi}{2} (y_{j}- y_{j-1})^2.
\end{equation}
This follows from the fact that
this holds
for $f_n$ approximating $\psi$ as in \cref{assu}.%
\kommentar{
it holds
\begin{equation*}
  f_n'(y_{j})( y_{j}- y_{j-1}) = f_n(y_{j}) - f_n(y_{j-1}) + \tfrac{1}{2}f_n''(s)(y_{j}-y_{j-1})^2 \geq
  f_n(y_{j}) - f_n(y_{j-1}) - \tfrac{C_\psi}{2} (y_{j}- y_{j-1})^2,
	\end{equation*}
	where $s$ is some appropriate intermediate point from the Taylor expansion
 and then using that $f_n\to\psi$ in $C^1_{\text{loc}}$.
}
 \kommentar{More in detail:
	By adding zero and rearranging we obtain
	\begin{align*}
&	\psi'(y_{j})( y_{j}- y_{j-1}) \geq 
  \psi(y_{j}) - \psi(y_{j-1}) - \tfrac{C_\psi}{2} (y_{j}- y_{j-1})^2 +\\
&\left[\left(\psi'(y_j)-f_n'(y_j)\right)(y_j-y_{j-1}) + (f_n(y_j)-\psi(y_j)) - (f_n(y_{j-1})-\psi(y_{j-1}))\right]
	\end{align*}
	and then taking the limit $n\to\infty$ which conserves the inequality.
	}
	\\
	Collecting terms and using the definition of the Ginzburg-Landau energy \cref{GinzburgLandauEnergyAniso} one finds
	\begin{equation}
          \left(\dfrac{1}{\tau_j} -\dfrac{C_\psi}{2
            }\right) \norm{y_{j}-y_{j-1}}^2 + 
          \left[\mathcal{E}(y_{j}) - \mathcal{E}(y_{j-1})\right] \leq 0
	\end{equation}
	and thus $\mathcal{E}(y_{j}) \leq \mathcal{E}(y_{j-1})$ if $\tau\leq 2/C_\psi$. 
\end{proof}
	The result from \cref{th:energy_stable_disc} can be applied to the discretizations assumed in \cref{th:existence_state_disc} and \cref{th:lipschitz_disc} since it provides the less restricting assumption on the step length $\tau$.
  
Finally let us comment on possible choices for the function $A$ for the quasilinear term and for the function $\psi$, in particular regarding the application to optimal control of anisotropic Allen-Cahn equations.
\begin{remark}\label{possiblepsi}~\\
    The \cref{assu} on $A$ are fulfilled,
  \revtwo{e.g., in the following cases:
    }
\begin{enumerate}
\item   \revtwo{
    $A\in C^1(\mathbb{R})$ is convex, absolutely $2$-homogeneous and satisfies $A(p)>0$ for $p\neq0$ as in \cite{Elliott1996}.       
Writing $A(p) =\tfrac{1}{2} \gamma(p)^2$, sufficient conditions on $\gamma $ for these properties can be found, e.g., in \cite{Graser2013}.
}

\revtwo{
  One can actually require only $A\in C^0(\mathbb{R})$. Then
the above results still hold for 
the solution of the variational inequality \eqref{proof_existence_state1} using the absolutely $2$-homogeneity in the estimate \eqref{proof_Elliot_1}.
    }
\kommentar{
* $ |A'(p)| \leq C|p|$ folgt, da $A'$ 1-homogen.\\
* aus Homogenität und $A(p)>0$ folgt $A(p)>c|p|^2$ mit $c>0$. Zusammen mit konvex folgt uniform convexity von $A$ also strong monotonicity von $A'$
}
	
\item
\revtwo{
$ A$ is given as   $A(p) =\tfrac{1}{2} \left( \sum_{l=1}^{L} {(p^T G_lp+\delta)}^{1/2}\right)^2
$ with symmetric positive definite  matrices $G_l\in \mathbb{R}^{d\times d}$ and $\delta \geq 0$, see \cite{BlMeNumeric} for details.
}

\revtwo{
If $\delta =0$, such anisotropies are studied in \cite{Barrett2014, bgn13} for the Allen-Cahn equation.
These may be regularized using $\delta >0$ to obtain $A \in C^2(\mathbb{R}^d)$ while loosing absolutely 2-homogeneity.
}
\end{enumerate}
  
The \cref{assu} on $\psi$ are fulfilled
in the following cases:
  \begin{enumerate}
  \item
$\psi\in C^2(\mathbb{R})$ 
is bounded from below, $\psi'' \geq -C_\psi$ for some $C_\psi \geq 0$ and $\lim_{t\to\pm\infty}\psi''(t) =+\infty$.\\
This can be shown by choosing a value
$x_{-1}>0$ large enough such that $\psi'' \geq 1, \psi' \geq 0$ on $[x_{-1},\infty)$. Then, with
	$x_n:= \text{argmin}_{x\in [x_{n-1}+1,\infty)}\psi''(x)$
	one can define the approximation on $[0,x_n]$ by $f_n:=\psi$, and for $x> x_n$ as $f_n(x):=\psi(x_n)+\psi'(x_n) (x-x_n) +\tfrac12 \psi''(x_n)(x-x_n)^2$. Respectively one can construct $f_n$ on $(-\infty, 0]$. 
\item
  $\psi$ is the double well potential $\psi(y) = \tfrac14(y^2-1)^2$, since then the conditions in 1. hold.
  
In addition
to the regularity $y \in L^\infty(0,T; H^1(\Omega))$
of the solution of~\cref{weakAC} that was shown in \cref{existence_state}, from the estimate \cref{eq:bounds_on_y_in_lemma_limit_from_discrete}
we also obtain
the regularity $y\in L^6(0,T; L^6(\Omega))$ for all space dimensions if we use this potential.
\item
$\psi$ is one of the following regularizations of the obstacle potential $\psi_{obst}$,
given by
$\tfrac12(1-x^2)  $ on $[-1,1]$ and $\infty$ elsewhere:
\begin{itemize}
	\item $\psi$ is the regularization considered in \cite{blowey_elliott_1991} for analyzing the solution of the isotropic Allen-Cahn or Cahn-Hilliard variational inequalities.
There  $\psi_{obst}$ is regularized to
$\psi \in C^2$ by a smooth continuation of $\psi_{obst}$ on $[-1,1]$ with a cubic polynomial in a neighborhood $\pm(1,1+\delta)$ 
and then with a quadratic polynomial
 (cf. formula (2.9) there).
\item $\psi$ is the Moreau-Yosida regularization of $\psi_{obst}$, i.e., $\psi \in C^1$ with 
$\psi(x)=\tfrac12(1-x^2) + s (\min\{x+1,0\})^2 + s (\max\{x-1,0\})^2 $ where the penalty parameter $s\in\mathbb{R}^+$ is possibly very large.
It is, e.g., used in \cite{HINTERMULLER2013810}
to study the optimal control of isotropic Allen-Cahn inequalities and to obtain a numerical approach.
\end{itemize}
\end{enumerate}
 \end{remark}

 \kommentar{   
It holds:\\
	Let $\psi\in C^2(\mathbb{R})$ with $\psi(x) \geq -\Lambda$ , $\psi''(x) \geq -\Lambda$ for all $x\in \mathbb{R}$ and $\Lambda\in \mathbb{R}$, $\lim_{x\to \pm \infty} \psi''(x) = \infty$. 
Then there exists a sequence $f_n$ that satisfies \cref{f_approximation_property}, i.e.
 $f_n \in C^2(\mathbb{R})$, 
 $-\Lambda \leq f_n \leq \psi$,
$f_n'' \geq -\Lambda$,
 $|f_n''| \leq C_n$,
$f_n \to \psi \quad \text{in } C^1_{\text{loc}}$.
}

\section{Existence of the optimal control in the time-discretized and in the continuous setting}
\label{sec:existence_control_problem}

Having shown the existence of solutions to the discretized \cref{scheme_state1} and time-continuous state equation \cref{weakAC} {that satisfy the bounds \cref{estim_for_y_tau} and \cref{eq:bounds_on_y_in_lemma_limit_from_discrete}, respectively}, we are able to develop the existence results of solutions to the pertinent control problems \cref{discrete_problem1} and \cref{eq1a}.

\begin{theorem}\label{excontroldisc}
  Let \cref{assu} be fulfilled
  and $\max_j \tau_j \coloneqq \tau < \tfrac1{(1+ 2C_\psi)}$ hold.
 Then for every $y_\Omega\in L^2(\Omega)$ the control problem \cref{discrete_problem1}--\cref{scheme_state1} has at least one solution in $U_\tau \times Y_\tau$.
\end{theorem}
\begin{proof}
  The requirements assure that \cref{th:existence_state_disc} is applicable and for every $u_\tau \in U_\tau$ we find a unique solution $S_\tau(u_\tau)=y_\tau\in Y_\tau$ of \cref{scheme_state1}. Since the feasible set $\{ (u_\tau, y_\tau)\mid y_\tau =S_\tau(u_\tau)  \mbox { for  } u_\tau \in U_\tau \}$
    is nonempty and the cost functional in \cref{discrete_problem1} is bounded from below we can deduce the existence of an infimum $\iota$ and of a minimizing sequence $((u_\tau^{(m)},y_\tau^{(m)}) )_m $
    with $\iota \coloneqq \lim_{m\to\infty} {J}(y_\tau^{(m)},u_\tau^{(m)})$.
    If $u_\tau^{(m)}\in L^2(0,T;L^2(\Omega))$ was unbounded so would be \revtwo{${J}(y_\tau^{(m)},u_\tau^{(m)})$} which would contradict its convergence to an infimum.
    Hence there exists a constant $\bar c_\tau > 0$ possibly depending on $\tau$ with $	\norm{u_\tau^{(m)}}_{L^2(0,T;L^2(\Omega))}\leq \bar c_\tau$ for all $m$
    and we can extract a weakly convergent subsequence denoted in the same way $ u_\tau^{(m)} \rightharpoonup u_\tau^\ast$ in $L^2(0,T;L^2(\Omega))$.
From \cref{th:existence_state_disc} we obtain independent from $m$
\begin{equation}
\label{eq:estimates_for_y_disc2}
\norm{y_\tau^{(m)}}_{L^\infty(0, T; H^1(\Omega))}
+ \norm{\psi'(y_\tau^{(m)})}_{L^2(0, T; L^2(\Omega))}
\leq C_{A,\psi,y_0}(\bar c_\tau) .
\end{equation}
\johannes{Given $y^{(m)}_\tau$, $y_\tau^\ast$ are determined by $N$ functions in $H^1(\Omega)$,}
this yields $y_\tau^{(m)} \to y_\tau^\ast$ in \revtwo{$L^\infty(0, T; L^2(\Omega))$ and in the pointwise sense},
 as well as
 \linebreak
 {$\psi'(y_\tau^{(m)} ) \allowbreak \rightharpoonup \psi'(y_\tau^\ast)$}
 in $L^2(0, T; L^2(\Omega))$ possibly for a subsequence.
Since $U_\tau$ is finite dimensional in time and due to the compact imbedding $L^2(\Omega)\hookrightarrow H^1(\Omega)'$ we obtain $u_\tau^{(m)} \to u_\tau^\ast$ in $L^2(0,T; H^1(\Omega)')$.
So the Lipschitz-continuity stated in \cref{th:lipschitz_disc} in addition yields
$y_\tau^{(m)} \to y_\tau^\ast $ in $L^2(0, T; H^1(\Omega))$.
Now we can take the limit in the state equation and obtain
	\begin{equation}
	\label{eq:limit_time_step_in_disc_existence_proof}
	(y_j^\ast - y_{j-1}^\ast, \varphi)
        + \tau_j(A'(\nabla y_j^\ast), \nabla \varphi)
        + \tau_j(\psi'(y_j^\ast), \varphi)
        = \tau_j (u_j^\ast, \varphi) \qquad j=1,\ldots,N.
	\end{equation}
        The convergence of the second term arises from the fact that $A': L^2(\Omega) \to L^2(\Omega)$ is a \revone{continuous} Nemytskii operator. From \cref{eq:limit_time_step_in_disc_existence_proof} we conclude that
        $y_\tau^\ast= S_\tau(u_\tau^\ast)$ and hence
        $(u_\tau^\ast,y_\tau^\ast)$ is feasible
        and its optimality follows by using the weak lower-semicontinuity of ${J}$.
\end{proof}
%
%
%
Similarly we can show the existence of the optimal control in the time continuous setting given the control-to-state operator $S:u\rightarrow y$ and the estimates
\cref{eq:bounds_on_y_in_lemma_limit_from_discrete}  for $y$  provided in the proof of \cref{existence_state}.

\begin{theorem}
	\label{existencecontrol}
 If \cref{assu} and $y_\Omega \in L^2(\Omega)$ hold, then
there exists a solution to the optimization problem \cref{eq1a}--\cref{weakAC}.
\end{theorem}
\begin{proof}
  \kommentar{
	From \cref{existence_state}, for each $u\in L^2(0,T;L^2(\Omega))$ we get the existence of a {unique} weak solution $y=S(u)$ in $L^2(0,T;H^1(\Omega))\cap H^1(0,T;L^2(\Omega))$ to \cref{weakAC}.
	Since $J$ is bounded from below, we deduce the existence of an infimum of $J$ over the feasible set of solutions $(y,u)$ of \cref{weakAC} and denote the minimizing sequence by $(y_m, u_m)_{m\in \mathbb{N}}$.
	The sequence $(u_m)_{m\in \mathbb{N}}$ is bounded
	in $L^2(Q)$ due to
	$	J(y_{m_k}, u_{m_k}) \geq \dfrac{\lambda}{2} \|u_{m_k}\|^2$  and by reflexivity of $L^2(Q) $ we can extract a weakly convergent subsequence with limit $\bar{u}$.	
	We now choose $\Lambda$ at least large enough to bound the sequence $(u_m)_{m\in \mathbb{N}}$. Then the estimates \cref{condlemma} are fulfilled independently of $m$.
	Consequently the estimate \cref{eq:bounds_on_y_in_lemma_limit_from_discrete} holds independently of $m$ for the solutions $y_m=S(u_m)$ of \cref{weakAC} constructed in the existence proof.
	Hence, due  to the uniqueness of the solutions we have for $y_m$
      }
As in the proof of \cref{excontroldisc} we obtain a minimizing sequence $(u_m,y_m)$ with $y_m=S(u_m)$ where $u_m$  is bounded and consequently providing a constant $\bar c$ such that \cref{eq:bounds_on_y_in_lemma_limit_from_discrete} holds independently of $m$, i.e.,
	\begin{equation}
	\norm{\partial_t y_m}_{L^2(0,T; L^2(\Omega))} + \norm{y_m}_{\revtwo{L^\infty(0,T; H^1(\Omega))}} + \norm{\psi'(y_m)}_{L^2(0,T; L^2(\Omega))} \leq \bar{c}.
	\end{equation}
        From this we get a subsequence $(u_m,y_m) $ with
        $u_m$ converging weakly to a $\bar{u}$ in
$L^2(0,T;L^2(\Omega))$, and
$y_m$ converging to a \revtwo{function $\bar{y} $ weakly-$\ast$ in
$L^\infty(0,T;H^1(\Omega))$, weakly in $H^1(0,T;H^1(\Omega)')$, and therefore using the lemma of Aubin-Lions strongly in
$ C([0,T]; L^2(\Omega))$ with $\bar{y}(0)=y_0$}
	and pointwise almost everywhere in $Q$.
	\kommentar{ due to $L^2(0,T;H^1(\Omega))\cap H^1(0,T;H^1(\Omega)') \hookrightarrow C([0,T]; L^2(\Omega))$.}%
	Moreover, $\partial_t y_m $ and $\psi'(y_m) $ converge \revtwo{to $\partial_t \bar{y}$ and $\psi'(\bar{y})$, respectively, in the weak topology of} $L^2(0, T; L^2(\Omega))$.
	\kommentar{ from boundedness in $L^2(0,T;L^2(\Omega))$.
		That the limit actually equals $\psi'(\bar{y})$ is due to the pointwise a.e. convergence $y_m \to \bar{y}$ (from strong convergence in $L^2(Q)$) and continuity of $\psi'$.}%
	In order to obtain $\bar y = S(\bar u)$ we need to be able to pass to the limit also in the $A'$-term of
	\cref{weakAC}.
	Given the fact that $A': L^2(Q) \to L^2(Q)$
	is a \revone{continuous} Nemytskii operator it is sufficient to show
	the strong convergence $\nabla y_m \to \nabla \bar{y}$ in $L^2(0,T;L^2(\Omega))$.
	Then finally, the weak lower-semicontinuity of $J$%
	\kommentar{(deduced from its continuity and convexity as a functional $L^2(0,T;H^1(\Omega))\cap H^1(0,T;L^2(\Omega)) \times L^2(Q) \to \mathbb{R}$, where in particular the continuous imbedding $L^2(0,T;H^1(\Omega))\cap H^1(0,T;H^1(\Omega)')\hookrightarrow C([0,T]; L^2(\Omega))$ is used)}
	provides $(\bar{y}, \bar{u})$ being a minimizer of $J$.
	
        The time derivative is monotone if $y_m(0)-\bar{y}(0) = 0$.%
        \kommentar {
 which follows from
	\begin{equation}
	\braket{ \partial_t y_m, y_m-\bar{y}} = \braket{ \partial_t y_m -  \partial_t \bar{y}, y_m-\bar{y}} + \braket{ \partial_t \bar{y}, y_m-\bar{y}} \geq \braket{ \partial_t \bar{y}, y_m-\bar{y}}.
      \end{equation}
      }
      Hence  we have
\linebreak
      $\braket{ \partial_t \bar{y}, y_m-\bar{y}}\leq\braket{ \partial_t y_m, y_m-\bar{y}} $,
      and $y_m=S(u_m)$ yields
	\begin{equation*}
	(A'(\nabla y_m), \nabla y_m - \nabla \bar{y})
	\leq 
	(u_m,y_m-\bar{y})
	- 
	(\psi'(y_m),y_m-\bar{y})
	- \braket{\partial_t \bar{y},y_m-\bar{y}}.
	\end{equation*}
	Recalling the convergence properties of $y_m$ together with $\|u_m\| + \|\psi'(y_m)\| \leq C$, the right hand side vanishes in the limit $m\to \infty$.	
	From strong monotonicity we obtain
	\begin{equation*}
	C\|\nabla y_m-\nabla \bar{y}\|^2
	\leq (A'(\nabla y_m), \nabla y_m - \nabla \bar{y}) - (A'(\nabla \bar{y}), \nabla y_m - \nabla \bar{y}),
	\end{equation*}
	where the second term on the right hand side vanishes in the limit by weak convergence and we have just shown that the limit of the first one can be bounded by $0$ from above. This finally yields the desired strong convergence of $\nabla y_m$ in $L^2(Q)$.
\end{proof}
Note that for the convergence $\nabla y_m \to \nabla y$ in $L^2(0,T;L^2(\Omega))$ we could not use \revone{the Lipschitz-estimate \cref{eq:Lipschitz_dependency_estimate} for the time continuous problem like we were able to use the analogous estimate \cref{eq:Lipschitz_discrete_inequ2} for the time discrete problem in} \cref{excontroldisc}. The reason is that in the time continuous case we do not have the analogon to the compact imbedding $L^2(\Omega)^N \hookrightarrow H^1(\Omega)'^N$ ($L^2(0,T;L^2(\Omega)) \hookrightarrow L^2(0,T;H^1(\Omega)')$ is not compact). Therefore we had to show the convergence more directly.\\
%
As in \cite{Herzog2011} for elliptic and in \cite{WachsmuthI} for parabolic problems we finally consider the convergence of the minimizers $u_\tau$ of the discretized problem to a minimizer of the problem in the continuous setting. Note that in \cite{WachsmuthI} a target function $y_Q$ is given over the whole time horizon or the solution has to have a higher regularity with respect to time. Having shown the strong convergence result $y_\tau(T,\cdot)\to y(T,\cdot)$ in \cref{existence_state} we still can derive the following result.
%
\begin{theorem}
	\revone{Let \cref{assu} be fulfilled and $y_\Omega \in L^2(\Omega)$ hold.}
	Consider a sequence of global optimal controls $(u_\tau, y_\tau)_\tau$ of \eqref{discrete_problem1} \revone{subject to} \eqref{scheme_state1} belonging to a sequence of discretizations with $\tau\to0$.
	Then there exists a subsequence with $u_\tau \to u$ in $L^2(0,T;L^2(\Omega))$ \revone{and with $y_\tau$ converging to $y=S(u)$ in the sense of \eqref{eq:various_convergences_of_y_tau}} where $(u, \revone{y})$ solves \eqref{eq1a} \revone{subject to}~\eqref{weakAC}.
\end{theorem}
\kommentar{die ganze Folge muß nicht konvergieren, da nicht unbedingt ein eindeutiges globales Minimum existiert.}
\begin{proof}
	First we choose an arbitrary $u^\ast \in L^2(0,T;L^2(\Omega))$ and a sequence $u^\ast_\tau\in U_\tau$ with $u^\ast_\tau \to u^\ast$ in $L^2(0,T;L^2(\Omega))$. 
	Hence $ y^\ast_{\tau}=S_\tau(u^\ast_\tau)$ is bounded in
	$ L^\infty(0, T; H^1(\Omega))$ due to \cref{estim_for_y_tau}.
	Now let $(u_\tau)_\tau$ be the sequence of global minimizers to \cref{discrete_problem1}
	subject to \cref{scheme_state1} \revone{and denote $y_\tau=S_\tau(u_\tau)$}.
	Then $J(y_\tau, u_\tau) \leq J(y^\ast_\tau, u^\ast_\tau) \leq c $ 
	implies that $(u_\tau)_\tau$ is bounded in $L^2(0,T;L^2(\Omega))$ and we deduce a subsequence with $u_\tau \rightharpoonup u$ in $L^2(0,T;L^2(\Omega))$.
	Then the bounds \cref{estim_for_y_tau} are fulfilled and \cref{existence_state} yields that
	we have the strong convergence $y_\tau(T,\cdot)\to y(T,\cdot)$ in $L^2(\Omega)$ \revone{where $y=S(u)$}.
	Respectively, given some arbitrary sequence $\tilde{u}_\tau$ with $\tilde{u}_\tau \to \tilde{u}$ in $L^2(0,T;L^2(\Omega))$ we obtain the latter also for \revone{$\tilde{y}_\tau=S_\tau(\tilde{u}_\tau)$} and \revone{$\tilde{y}=S(\tilde{u})$}.
	This yields
	\begin{equation}
		\label{eq:conv_disc_to_cont_inequality_chain}
		J(y,u) \leq \liminf_{\tau\to 0} J(y_\tau, u_\tau) \leq \revtwo{\limsup_{\tau\to 0} J(y_\tau, u_\tau) \leq \lim_{\tau\to 0} J(\tilde{y}_\tau, \tilde{u}_\tau)} = J(\tilde{y}, \tilde{u}).
	\end{equation}
	Since $\tilde{u}$ was arbitrary this yields the global optimality of $u$. Plugging in $\tilde{u} = u$ yields the convergence $\|u_\tau\| \to \|u\| $
	and therefore
	with the weak convergence also
	the strong convergence $u_\tau \to u$ in $L^2(0,T;L^2(\Omega))$.
\end{proof}
\begin{remark}
If instead of the cost functionals \cref{eq1a} and \cref{discrete_problem1} one considers the cost functionals with a  target function $y_Q$ in $L^2(0,T;L^2(\Omega))$ given over the whole time horizon
	\begin{equation}
	\label{eq:distributed_cost_cont}
	J(y,u) \coloneqq \dfrac{1}{2} \norm{y - y_Q}^2_{L^2(Q)} + \dfrac{\lambda}{2} \norm{u}^2_{L^2(Q)}
	\end{equation}
	and its discrete counterpart
	\begin{equation}
	\label{eq:distributed_cost_disc}
	{J}_\tau({y_\tau},{u_\tau}) \coloneqq \dfrac{1}{2}\sum_{j=1}^N \tau_j \|y_j-y_{Q,j}\|^2 + \dfrac{\lambda}{2}
	\sum_{j=1}^N \tau_j \|u_{j}\|^2,
	\end{equation}
	with $y_{Q,\tau}\in Y_\tau$ and $y_{Q,\tau}\to y_Q$ in $L^2(0,T;L^2(\Omega))$, the theorems of this section still hold true with proofs following the same lines.
\end{remark}

\section*{Acknowledgements}
The authors gratefully acknowledge the support by the RTG 2339 “Interfaces, Complex Structures, and Singular Limits” of the German Science Foundation (DFG).

\printbibliography
\clearpage
\end{document}